\documentclass[12pt]{amsart}
\usepackage{amsmath,amssymb}
\usepackage{graphicx,epsfig,subfigure,psfrag}
\begin{document}

%
%
\newtheorem{theorem}{Theorem}
\newtheorem{proposition}[theorem]{Proposition}
\newtheorem{lemma}[theorem]{Lemma}
\newtheorem{corollary}[theorem]{Corollary}
\newtheorem{definition}[theorem]{Definition}
\newtheorem{remark}[theorem]{Remark}
\numberwithin{equation}{section}
\numberwithin{theorem}{section}
\newcommand{\be}{\begin{equation}}
\newcommand{\ee}{\end{equation}}
\newcommand{\tex}{\textstyle}
\newcommand{\re}{\mathbb{R}}
\newcommand{\n}{\nabla}
\newcommand{\ren}{\mathbb{R}^N}
\newcommand{\iy}{\infty}
\newcommand{\pa}{\partial}
\newcommand{\ms}{\medskip\vskip-.1cm}
\newcommand{\mpb}{\medskip}
\newcommand{\inA}{\quad \mbox{in} \quad \ren \times \re_+}
\newcommand{\inB}{\quad \mbox{in} \quad}
\newcommand{\inC}{\quad \mbox{in} \quad \re \times \re_+}
\newcommand{\inD}{\quad \mbox{in} \quad \re}
\newcommand{\forA}{\quad \mbox{for} \quad}
\newcommand{\whereA}{,\quad \mbox{where} \quad}
\newcommand{\asA}{\quad \mbox{as} \quad}
\newcommand{\andA}{\quad \mbox{and} \quad}
\newcommand{\withA}{,\quad \mbox{with} \quad}
\newcommand{\orA}{,\quad \mbox{or} \quad}
\newcommand{\atA}{\quad \mbox{at} \quad}
\newcommand{\onA}{\quad \mbox{on} \quad}
\newcommand{\ef}{\eqref}
\newcommand{\mc}{\mathcal}
\newcommand{\mf}{\mathfrak}
\newcommand{\ssk}{\smallskip}
\newcommand{\LongA}{\quad \Longrightarrow \quad}

\newcommand{\BB}{{\bf B}}
\newcommand{\Am}{{\bf A}_{2m}}
\newcommand{\bL}{\BB^*}
\newcommand{\bLs}{\BB}
\renewcommand{\a}{\alpha}
\renewcommand{\b}{\beta}
\newcommand{\g}{\gamma}
\newcommand{\ka}{\kappa}
\newcommand{\G}{\Gamma}
\renewcommand{\d}{\delta}
\newcommand{\D}{\Delta}
\newcommand{\e}{\varepsilon}
\renewcommand{\l}{\lambda}
\renewcommand{\o}{\omega}
\renewcommand{\O}{\Omega}
\newcommand{\s}{\sigma}
\renewcommand{\t}{\tau}
\renewcommand{\th}{\theta}
\newcommand{\z}{\zeta}
\newcommand{\wx}{\widetilde x}
\newcommand{\wt}{\widetilde t}
\newcommand{\noi}{\noindent}
\newcommand{\lb}{\left (}
\newcommand{\rb}{\right )}
\newcommand{\lsb}{\left [}
\newcommand{\rsb}{\right ]}
\newcommand{\lab}{\left \langle}
\newcommand{\rab}{\right \rangle }
\newcommand{\gap}{\vskip .5cm}
\newcommand{\bz}{\bar{z}}
\newcommand{\bg}{\bar{g}}
\newcommand{\Ba}{\bar{a}}
\newcommand{\bt}{\bar{\th}}
\def\com#1{\fbox{\parbox{6in}{\texttt{#1}}}}

\title{\bf On oscillations of solutions of the fourth-order thin film equation
near heteroclinic bifurcation point}

\author
{V.A. Galaktionov}

\address{
 Department of Math. Sci., University of Bath,
 Bath, BA2 7AY}
\email{masvg@bath.ac.uk}


\date{\today}



\begin{abstract}

The behaviour of   solutions of the Cauchy problem for the 1D
fourth-order thin film equation (the TFE--4)
\[
u_t = - (|u|^n u_{xxx})_{x} \quad \mbox{in} \quad \re \times
\re_+,
\]
where, mainly, $n \in (\frac 32,2]$, is studied.  Using the
standard mass-preserving (ZKB-type) similarity solution of the
TFE--4, a third-order ODE for the profile $f(y)$ occurs:
\[ 
  \begin{matrix}
u_{\rm s}(x,t)=t^{- \frac 1{4+n}}f(y), \,\,\, y= \frac
x{t^{1/(4+n)}} \LongA -(|f|^n f''')' + \frac 1{4+n}\, (y f)'=0
\ssk \\
 \LongA -|f|^n f'''+ \frac 1{n+4}\, y f=0.
 \end{matrix}
 \] 
 For the Cauchy problem, the oscillatory behaviour
of $f(y)$ near the interface was studied in \cite{Gl4, PetI}, etc.
It was shown that a periodic oscillatory component $\varphi(s)$,
describing the actual sign-changing behaviour of $u_{\rm s}(x,t)$
close to finite interfaces, exists up to a critical {\em
homoclinic bifurcation} exponent $n_{\rm h}$, i.e., for
  \[ 
   \mbox{$
  0<n<n_{\rm h} \in (\frac 32,n_{\rm +}), \quad \mbox{where}
   \quad n_{\rm +}= \frac 9{3+\sqrt 3}=1.9019238... \, .
    $}
  \] 
Careful numerics show that  $n_{\rm h}= 1.75987... \, .$

In the present paper, a non-oscillatory behaviour of $f(y)$ for $n
\in [n_{\rm h},2]$, i.e., above the heteroclinic value, is
revealed by a combination of analytical and numerical methods. In
particular, this implies that, for $n \in [n_{\rm h},2]$ (and, at
least, up to $3$), the Cauchy setting coincides with the standard
zero-height-angle-flux free-boundary one (with non-negative
solutions), studied in detail in many well-known papers since the
end of 1980s, including Bernis--Friedman's seminal paper
\cite{BF1} of 1990.

\end{abstract}

\maketitle

\section{\sc Introduction: TFE--4, first similarity solution, and main result}
\label{Sect1}


\subsection{The TFE--4 and first mass-preserving similarity solution}

We study the first mass-preserving {\em self-similar solution}
 of the fourth-order {\em thin film equation}
 (the TFE--4) in one dimension
  \be
 \label{GPP}
 u_t = - (|u|^n u_{xxx})_x \quad \mbox{in} \,\,\, \re \times \re_+,
   \ee
 where, in the most of the cases,  $n \in (\frac 32,2]$ is a fixed exponent.
 Equation \eqref{GPP} is written for solutions of changing sign, which is an essential
 feature
for the Cauchy problem (CP) to be considered, though, in a
supercritical  parameter $n$-range (see below), solutions turn out
to be non-negative, as happens for the standard and much more
well-studied (at least, since 1980s) free-boundary problem (FBP)
with zero-height, zero-contact angle, and zero-flux conditions;
see references in \cite{Gl4} and in several other papers mentioned
below.


 Thus, we will need well-known similarity solutions of the
  TFE--4 \ef{GPP}.
 Source-type mass-preserved (i.e., of a ZKB-type)
similarity solutions of the FBP for (\ref{GPP}) for arbitrary $n$
were studied in \cite{BPelW92} for $N=1$ and in \cite{BFer97} for
the equation in $\ren$. More information on similarity and other
solutions can be found in \cite{Beck05, BerHK00, BerHQ00, Bow01,
   Car07, CarrT02}.
 TFEs admit
non-negative solutions constructed by ``singular" parabolic
approximations of the degenerate nonlinear coefficients. We refer
to the pioneering paper by Bernis and Friedman \cite{BF1} (1990)
and to various extensions in \cite{Ell96, EllS, Gia08, Govor05,
Green78, LPugh, Oron97, WitBerBer} and the references therein.
 See also
\cite{Gr95} for mathematics of solutions of the FBP and CP of
changing sign (such a class of solutions of the CP will be
considered later on).

In both the FBP and the CP, the source-type solutions of the
TFE--4 (\ref{GPP}) are
 \be
\label{uss1}
 \tex{
 u_{\rm s}(x,t) = t^{- \frac 1{4+n}} f(y), \quad y=
 \frac x {t^{1/(4+n)}},
 }
  \ee
 where, on substitution of \ef{uss1} into \ef{GPP}, $f(y)$ solves a divergent fourth-order
 ordinary differential equation (ODE):
 \be
\label{Od11}
 \tex{
 -(|f|^n f''')' + \frac 1{4+n} \, (y f)' =0.
 }
  \ee
On integration once, by using the continuity of the flux function
$|f|^n f'''$ at the interface, where $f=0$, one obtains a simpler
third-order ODE of the form
 \be
 \label{ODE3}
  \tex{
  |f|^n f'''= \frac 1{4+n} \, y f.
   }
   \ee
 In view of existence of evident scaling, we, first, get rid of the
 multiplier $\frac 1{4+n}$, and, next, pose the symmetry condition
 at the origin and the last normalizing one to get
  \be
  \label{up2}
|f|^n f'''=  y f \quad \mbox{for} \quad y>0, \quad \mbox{with}
\quad f'(0)=0, \andA f(0)=1.
 \ee
 Later on, for asymptotic expansions close to the interface point,
  assuming that such a solution has its interface as some finite
  $y_0>0$ (so that $f(y) \equiv 0$ for $y \ge y_0$),
  we perform a reflection,
   \be
   \label{up3}
   y-y_0 \mapsto y \LongA |f|^n f'''= - (y_0-y) f \quad \mbox{for small} \quad
 y>0.
  \ee

From now on, we will mainly concentrate on solutions and the
corresponding profiles $f(y)$ for the Cauchy problem only, though,
for some  supercritical range from $n \in \big(\frac 32,2\big]$
and even further, as will turn out, we  will conclude that this
covers the FBP setting as well.

For the Cauchy problem, it was formally  shown
 that there exists an
oscillatory similarity profile of (\ref{Od11}), that are
infinitely oscillatory as $y \to y_0^-$, for not that large
exponents\footnote{Not surprisingly and obviously, this range
includes $n=0$, i.e., classic analytic solutions of the CP of the
{\em bi-harmonic} (parabolic) equation $u_t=-u_{xxxx}$ obtained
from \ef{GPP}  by formally passing to the limit $n \to 0^+$. As we
have suggested earlier, this natural fact can be used as a
``proper" definition of CP-solutions of \ef{GPP}: these are those
that can be continuously deformed to the bi-harmonic solutions of
the CP using this homotopic path. For FBP-solutions, this is not
possible \cite{BerHK00}: non-negative FBP-solutions for $n>0$
converge to similar FBP-ones for $n=0$, which are not
oscillatory.}
 \be
 \label{up0}
n \in (0,n_{\rm h}) \whereA n_{\rm h}=1.75987...;
 \ee
 see \cite{Gl4}.
For $n \in (0,1)$, existence and uniqueness (up to the
mass-scaling) of an oscillatory $f(y)$  are straightforward
consequences of the result in \cite{BMcL91} on oscillatory
solutions for the fully divergent fourth-order {\em porous medium
equation}; see \ef{mm1} below. See also some details in
\cite[\S~3.7]{GSVR} and interesting oscillatory features of
similarity dipole solutions of the TFE (\ref{GPP}) observed in
\cite{BW06}. Stability of such sign changing similarity solutions
is quite plausible but was not proved rigorously being an open
problem.

Thus, in the present paper, we concentrate on the still unknown
parameter range
 \be
 \label{up1}
 n \in [n_{\rm h},2].
 \ee

 Our main goal is to show that, in this parameter range, the
 above similarity
 solutions, as solutions  of the Cauchy problem, {\em are not oscillatory}, and,
 moreover, do not change sign in small neingbourhoods
 of finite interfaces. In addition, we then claim that, for such
 $n$'s in \ef{up1}, {\em similarity} (and, most plausibly, not only those)
{\em solutions of the CP and the FBP} (see a number of papers
mentioned above) {\em coincide}.

\section{Heteroclinic subcritical and supecritical ranges}
 \label{SectHet}

We begin with a brief explanation of known results on what happens
in the heteroclinic subcritical range \ef{up0}, which is important
for our final conclusions.

\subsection{Local oscillatory structure near interfaces for $n \in (0, n_{\rm h})$}

 Namely,
 it is known
 that the asymptotic behaviour of $f(y)$ satisfying equation (\ref{ODE3})
near the interface point as $y \to y_0^- >0$ is given by the
expansion
 \be
\label{LC11}
 \mbox{$ f(y) = (y_0-y)^{\mu} \varphi(s) , \quad s
=\ln(y_0-y), \quad \mu = \frac 3 n,
 $}
   \ee
where, on substitution to \ef{up2}, after scaling, the {\em
oscillatory component}
 $\varphi$ 
 satisfies the following autonomous ODE,
  where an exponentially small as $s \to -\infty$
 term (occurring by setting $y=y_0 - {\mathrm e}^s$) is omitted:
 \be
\label{m=2.11} 
 \textstyle{\varphi''' + 3(\mu-1) \varphi''  + (3
\mu^2 - 6 \mu +2) \varphi'+ \mu(\mu-1)(\mu-2)\varphi +  \frac
\varphi {|\varphi|^{n}}=0.}
  \ee
It was shown that here exists a stable (as $s \to +\iy$)  changing
sign periodic solution $\varphi_*(s)$ of (\ref{m=2.11}). The
bifurcation value $n_{\rm h}$ was obtained in \cite{Gl4} by some
analytic and numerical evidence showing that, as $n \to n_{\rm
h}^-$, the ODE (\ref{m=2.11}) exhibits a usual heteroclinic
bifurcation, where a periodic solution is generated from a
heteroclinic path of two constant equilibria, \cite[Ch.~4]{Perko}.
 According to (\ref{LC11}), this
gives similarity profiles of changing sign, which being extended
by $f(y) \equiv 0$ for $y>y_0$ form a
 compactly supported solution\footnote{As usual, such a regularity is based of the actual smoothness of the algebraic envelope
 $(y_0-y)^\mu$
  in
\ef{LC11} and does not involve ``transversal zeros" nearby, at
which such a smoothness  breaks down.}
   $f \in C^\a$ in a neighbourhood of $y=y_0$,  with $\a \sim \frac
 3n$. Notice that $\a \to +\infty$ as $n \to 0^+$, so the
 regularity at $y=y_0$ improves to $C^\infty$ at $n=0^+$ (and eventually, to an analytic profile $F(y)$
 for $n=0$ for all $y>0$, being the ``Gaussian" kernel of the fundamental
 solution $b(x,t)=t^{-\frac 14} F(x/t^{\frac 14})$ of the linear parabolic operator $D_t+D_x^4$).

\ssk

{\bf Remark: TFE in $\ren$.} The above conclusions, at least,
formally can be applied to the Cauchy problem for the TFE--4  in
$\ren$:
 \be
 \label{TFERN}
 u_t=- \n \cdot (|u|^n \n \D u) \quad \mbox{in} \quad \ren \times
 \re_+,
  \ee
  with a bounded compactly supported data $u_0(x)$. Namely, at any
  point of a sufficiently smooth interface surface, at which there
  exists a unique tangent hyperplane $\pi$, the same 1D oscillatory
  behaviour is expected to occur in the direction of the inward
  unit normal to $\pi$. In other words, the principal phenomenon of
  existence of oscillatory sign-changing behaviour near the
  smooth interfaces is essentially 1D, excluding some special and
  ``singular" non-generic cases.

  On the other hand, at possible points of singular ``cusps" at
  the interfaces (anyway, expected to be non-generic), such an
  easy approach is not applicable, though almost nothing is known
  about those singular phenomena for \ef{TFERN} for any $N \ge 2$.

  \ssk

Note that the first  results on the oscillatory behaviour of
similarity profiles for fourth-order ODEs related to the
source-type solutions of the divergent parabolic PDE of the porous
medium type (the PME--4)
  \be
  \label{mm1}
  u_t = - (|u|^{m-1} u)_{xxxx} \quad (m>1)
   \ee
 were obtained in \cite{BMcL91}.
 It turns out that these results can be applied to the rescaled TFE
 (\ref{Od11}), but for $n \in (0,1)$ only. As shown in \cite{Gl4}, the corresponding ODEs
 for source-type similarity solutions for
 (\ref{GPP}) and (\ref{mm1}) then coincide after a change of parameters $n$ and
 $m$:
  $$
  \mbox{$
   m= \frac 1{1-n}, \quad n \in (0,1).
    $}
    $$
 Some results on existence and multiplicity of
periodic solutions of ODEs such as (\ref{m=2.11}) are known
\cite{Gl4}, and often lead to a number of open problems.
Therefore,  numerical and some analytic evidence remain key,
especially for sixth and higher-order TFEs; see \cite{GBl6} and
\cite[Ch.~3]{GSVR}.

Returning to the behaviour as $s \to -\iy$ for \ef{m=2.11}, i.e.,
approaching the interface, we express the above results as
follows: in the heteroclinic subcritical range, there exists a 2D
bundle of asymptotic solutions of (\ref{m=2.11}) in $\re$ having
the expansion (\ref{LC11}) (an exponentially small term is again
omitted)
 \be
 \label{2par}
 f(y)=(y_0-y)^{\frac 3n} [\varphi_*(s+s_0)+...],
  \quad \mbox{with parameters $y_0>0$ and  $s_0 \in \re$},
  \ee
  where we take into account the phase shift $s_0$ of the periodic orbit
 $\varphi_*(s)$. Here, we  should treat $y_0$ as a free parameter, though, in the above setting, it is
 fixed
 by conditions at $y=0$. Recall that, as $s \to -\iy$, the
 periodic solution $\varphi_*(s)$ is thus unstable by the obvious
 reason: a 1D unstable manifold is then generated by small
 perturbations of $y_0$.

Therefore, matching the 2D bundle (\ref{2par}) with two symmetry
conditions at the origin in (\ref{up2}) leads to a formally
well-posed problem of 2D matching.
However, using this matching approach directly is a difficult open
problem (for $n \in (1, n_{\rm h})$, where existence of $f(y)$ was
not properly proved; recall that, for $n=1$, there is a pure
algebraic approach to construct
 such a unique oscillatory profile \cite{Gl4}).
 Note that, in
the case of an analytic dependence on parameters involved, such a
problem cannot possess more than a countable set of isolated
solutions.

 \subsection{$n \in (\frac 32,2]$: numerical results}

 Before deriving proper asymptotic expansions in the
 supercritical range in the next section, we present some numerical evidence that, at
 $n = n_{\rm h}$, oscillations of $f(y)$, as a solution of \ef{up2}, near the interface disappear.

As usual, for numerics, we use the regularization in the
third-order operator in (\ref{up2}) by replacing
 \be
 \label{Reg1}
  |f|^n \mapsto (\e^2 + f^2)^{\frac n2}, \quad \mbox{with} \quad
 \e= 10^{-11},
   \ee
 which is sufficient for the accuracy required; see below.
 Note that, for our purposes in the Cauchy problem, we principally cannot use the
 pioneering  FBP-regularization, introduced in Bernis--Friedman
 \cite{BF1} (not analytic as in (\ref{Reg1}), except $n=2$, to be treated
 specially),
 leading, in many cases, to non-negative solutions.
 However, we expect that, for $n \in [n_{\rm h},2]$, both analytic
 and non-analytic $\e$-regularizations lead to the same solutions,
 i.e., the Cauchy and FBP settings coincide (though proving that
 rigorously is very difficult).

In our numerics by the {\tt MatLab}, we were often obliged to keep
a high accuracy, with both tolerancies $\sim 10^{-12}$, and the
total number of points on the interval $y \in (0,3)$ of about
50000--100000, though using the simplest ODE solver for the Cauchy
problem {\tt  ode45} took, sometimes, 5--10 minutes for each
shooting from $y=0$ for each fixed value of $n$.

There appears a simple one-parameter shooting problem for the ODE
\ef{up2} from $y=0$ with the only parameter
 \be
 \label{Reg55}
 f''(0)= \mu \in \re.
 \ee
 Then, varying $\mu<0$, we are trying to approach as close as
 possible to the corresponding finite interface at some $y_0>0$,
 and, eventually, to see whether $f(y)$ changes sign or not in a
 small vicinity of $y_0^-$. A few of such results are presented
 below.

\begin{figure}
\centering
\includegraphics[scale=0.65]{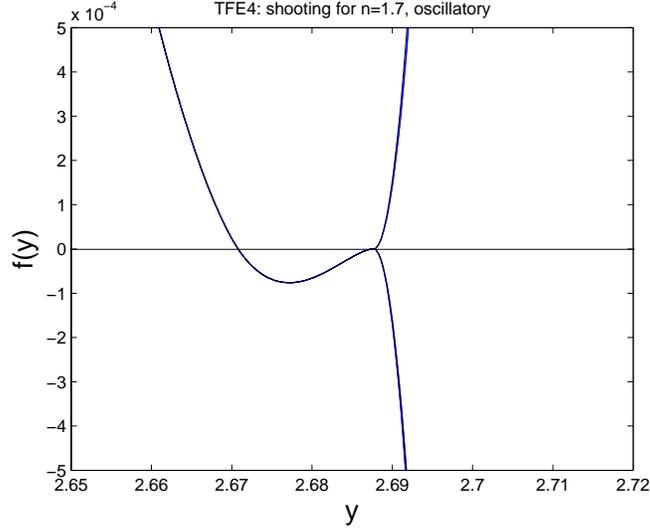}                     
 \vskip -.4cm
\caption{\rm\small The first zero, close to $y=2.67$, of $f(y)$
near $y=y_0$ for $n=1.7$.}
 \label{F1}
\end{figure}

Thus, let us begin with the oscillatory range \ef{up0}. In Figure
\ref{F1}, our numerics catch a changing sign behaviour as $y \to
y_0^-$ for $n=1.7< n_{\rm h}$. Our shooting from $y=0$, with
$f(0)=1$, even with the above mentioned accuracy and with about
$10^5$ points, allows us to see just a single zero near $y_0^-$.
We do not think that, in a few minutes of PC time, as here, the
next, to say nothing of other zeros of $f(y)$, near $y_0^-$ can be
viewed by such a shooting.

In the next Figure \ref{F2}, we start to approach the heteroclinic
bifurcation value $n_{\rm h}$ from below. Namely, we take
$n=1.75$, which is still slightly  less. We again were able to see
the first zero, though to reach that, we increased the number of
points, but still observed a ``non-smooth" profile $f(y)$ close to
the interface. We see that a negative hump in this figure, which
is already of the order $\sim 10^{-7}$, so one can expect that
further approaching $n_{\rm h}$ from below will require a
completely different and more powerful  computer tools.

 \begin{figure}
\centering
\includegraphics[scale=0.65]{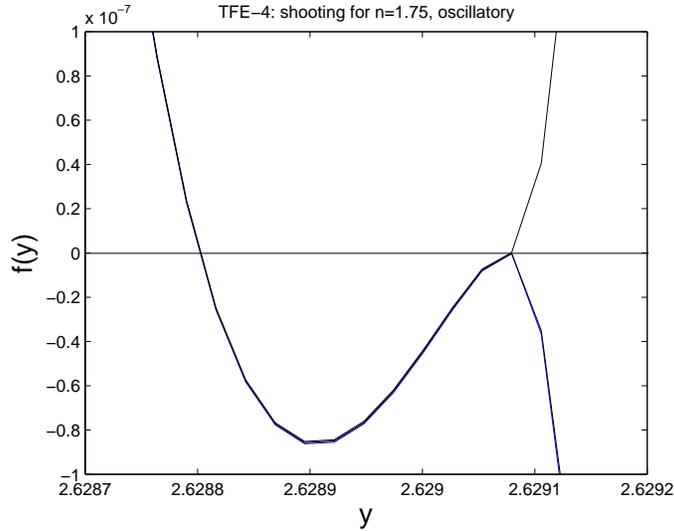}                     
 \vskip -.4cm
\caption{\rm\small The first zero, close to $y=2.6288$, of $f(y)$
near $y=y_0$ for $n=1.75 < n_{\rm h}$. Here, for the last profile,
$f''(0)=-0.434097009$.}
 \label{F2}
\end{figure}

\begin{figure}
\centering
\includegraphics[scale=0.65]{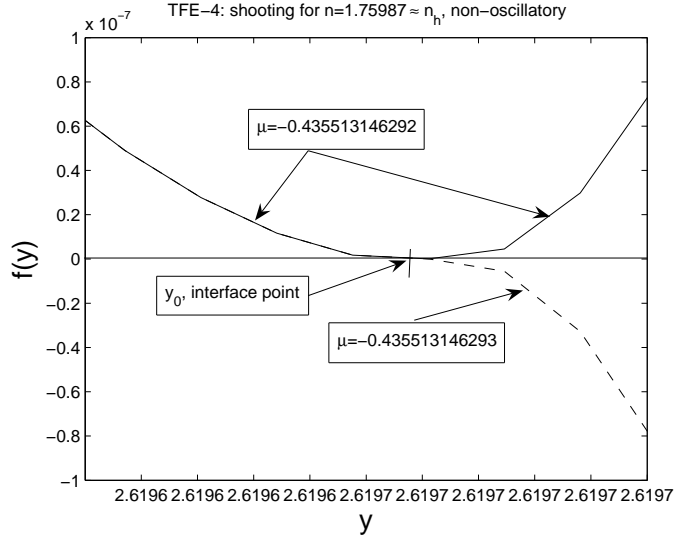}                     
 \vskip -.4cm
\caption{\rm\small Non-oscillatory $f(y)$ near $y=y_0$ for
$n=1.75987 \approx n_{\rm h}^- $. Here, for the last profile,
$\mu= f''(0)=-0.435513146293$.}
 \label{Fnh}
\end{figure}

We next take $n$ very close to $n_{\rm h}$, where
 \be
 \label{nnn1}
 n=1.75987, \quad f''(0)=\mu=-0.435513146293, \quad y_0=2.6197.
 \ee
 In this case we enhance our numerics by taking 200000 points on
 the interval
 $(0,2.7)$ to get the results in Figure \ref{Fnh}. Thus, up to the
 accuracy $10^{-7}$ (and more, as easily seen), $f(y)$ {\em is not
 oscillatory} for $y \approx y_0^-$, as our analytical predictions
 said earlier.

\begin{figure}
\centering
\includegraphics[scale=0.65]{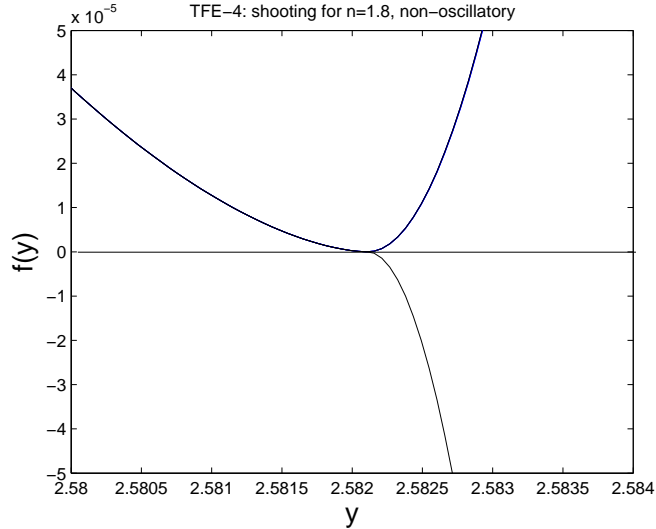}                     
 \vskip -.4cm
\caption{\rm\small Non-oscillatory  $f(y)$ for $n=1.8 > n_{\rm h}
$.}
 \label{F3}
\end{figure}

\begin{figure}
\centering
\includegraphics[scale=0.65]{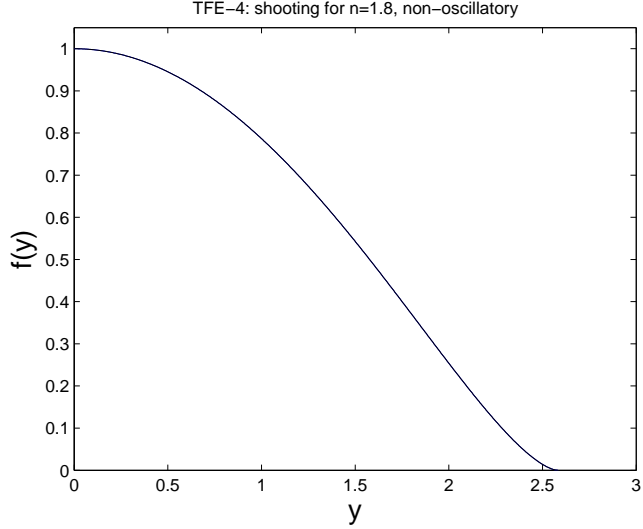}                     
 \vskip -.4cm
\caption{\rm\small The non-oscillatory profile $f(y)$ for $n=1.8 >
n_{\rm h} $.}
 \label{F4}
\end{figure}

Now, consider the supercritical range above $n_{\rm h}$.  Figure
\ref{F3}, constructed for $n=1.8> n_{\rm h}$ (but rather close to
it), shows that, with the prescribed accuracy, we do not see any
sign  changes near the interfaces, unlike the previous figures. In
Figure \ref{F4}, we show a global structure of this $f(y)$, which
looks like being non-negative, but, to verify that, one needs more
``microscopic" structure near the interface given in the previous
figure.

\begin{figure}
\centering
\includegraphics[scale=0.65]{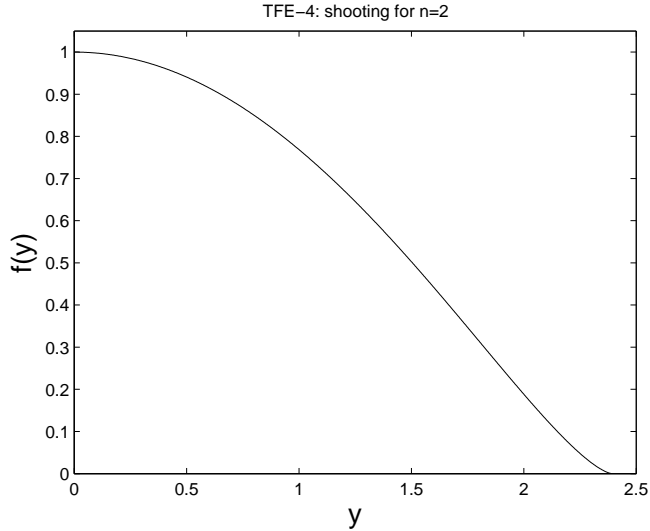}                     
 \vskip -.4cm
\caption{\rm\small The non-negative profile $f(y)$ for $n=2$.}
 \label{F5}
\end{figure}

Finally, for completeness, in Figure \ref{F5}, we present  a full
view of $f(y)$ for the special case $n=2$, where both analytic and
non-analytic (for the FBP) regularizations coincide. Again,
``microscopically", no sign changes of $f(y)$ were observed.

We do not expect (and did not see) any essential changes by
further increasing $n \in (2,3)$ until the special $n=3$ to be
treated next.

\subsection{$n=3$: logarithmically perturbed linear asymptotics
and existence of $f(y)$ for the CP}

\begin{figure}
\centering
\includegraphics[scale=0.65]{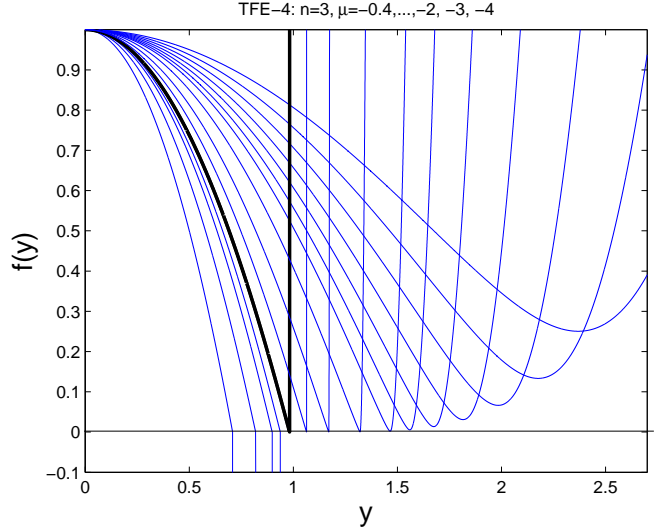}                     
 \vskip -.4cm
\caption{\rm\small Non-oscillatory $f(y)$  for $n=3$, the
bold-face line. For it, $\mu=f''(0)=-2.274...$\, ,
$y_0=0.943...$\, .}
 \label{Fn3}
\end{figure}

For $n=3$, similar to the case $n= \frac 32$ for the
FBP\footnote{Recall that, for the CP, the case $n=\frac 32$ is not
special, and for all $n \in (0,n_{\rm h})$ the solutions $f(y)$
are equally oscillatory near the interface, governed by a periodic
oscillatory component.},
 we observe a logarithmically perturbed linear asymptotics (we do not reveal the second term, since $n=3$ is currently
 out of our main interests): as $y \to y_0^-$,
  \be
  \label{log1}
   \tex{
  f(y)= \frac 3{\sqrt 2}\, (y_0-y) |\ln(y_0-y)|^{\frac 13}+... \,.
   }
   \ee
 The corresponding shooting of $f(y)$ is presented in Figure
 \ref{Fn3}. Obviously, such similarity solutions with the expansion \ef{log1} do
 not satisfy the zero contact angle condition (however, other FBPs can be posed), but can be proper
 solutions of the Cauchy problem, so that this behaviour at
 interfaces can be that of a {\em maximal regularity}.

\subsection{$n=4$: nonexistence of $f(y)$}

\begin{figure}
\centering
\includegraphics[scale=0.65]{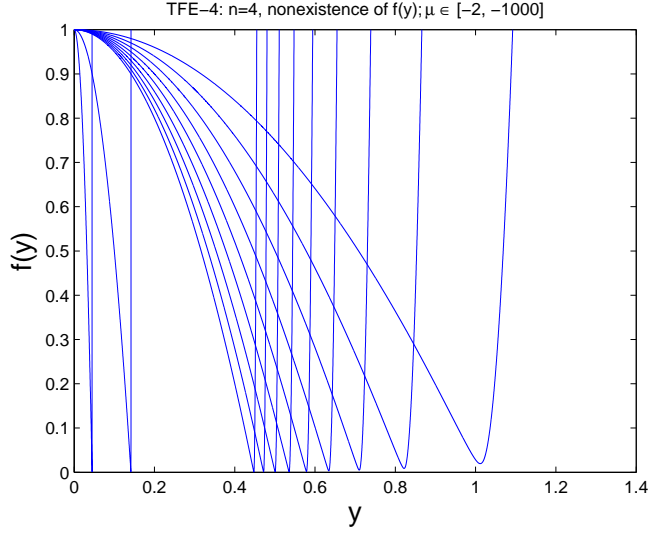}                     
 \vskip -.4cm
\caption{\rm\small Nonexistence of  $f(y)$ with a finite interface
for $n=4$.}
 \label{Fn4}
\end{figure}

As a simple illustration of the general results for the TFE--4 in
\cite{BF1}, saying that solutions of this PDE are strictly
positive for all $n$ large enough, consider our ODE \ef{up2} for
the next integer $n=4$. The results of shooting from $y=0$ are
shown in Figure \ref{Fn4}: no solutions for $\mu=f''(0)$ from $-2$
to $-1000$ can reach the zero level $f=0$. It is  clear that one
does not need any proof (though it is already available in
\cite{BF1}, and in much more general PDE fashion).

Anyway, as an ODE illustration again, we see from \ef{up3},
$y_0=1$ by scaling, multiplying by $f'$ and integrating that,
keeping the main singular terms,
 \be
 \label{JJJ1}
  \tex{
  f'f'' = \int (f'')^2 + \frac 12 \, \frac 1{f^2} + ... \ge \frac 12 \, \frac 1{f^2} \, .
   }
   \ee
However, the resulting ``majorising" ODE
 \be
  \label{JJJ2}
   \tex{
   f'f'' = \frac 12 \, \frac 1{f^2}
   }
   \ee
   does not admit increasing positive solutions for $y>0$ small (here, $y_0=0$ and $f(0)=0$, recall).
   Indeed, multiplying by $f'$ and integrating yields:
    \be
    \label{JJJ4}
     \tex{
      \frac 13 \, (f')^3= - \frac 12 \, \frac 1 f+... <0 \quad
      \mbox{for any small} \quad f>0.
      }
      \ee

\section{Positive expansions near interface for $n \in (\frac 32,3)$}
 \label{Snon}

Thus, we consider the ODE \ef{up3} close to the interface $y_0=1$,
by scaling. Assuming that $f(y)>0$ for small $y>0$ yields the
following equation:
 \be
 \label{E1}
 f^{n-1}f'''=-1+y.
  \ee
  First of all, we cannot use the standard ``parabolic" expansion from dozens of recent papers,
  since it is suitable for the FBP for $n \in (0, \frac 32)$.

Therefore, we are looking for an expansion about a different
profile (an explicit solution of \ef{E1} with $-1$ on the RHS):
 \be
 \label{E2}
  \tex{
 f_0(y)= B_0 y^m \whereA m= \frac 3n \quad \mbox{and}
 \quad
 B_0^n=\big[ \frac 3n \big(\frac 3n-1 \big)\big(2- \frac 3n\big)
 \big]^{-1},
  }
  \ee
which makes sense for any
 \be
 \label{E3}
  \tex{
  n \in \big( \frac 32,3 \big).
  }
  \ee

  We next use an algebraic perturbation of \ef{E2} by setting
   \be
   \label{E4}
    \tex{
    f(y)= B_0 y^m+D y^l+..., \quad \mbox{where} \quad l>  m=
    \frac 3n
    }
    \ee
    and $D \in \re$ is an arbitrary parameter (recall that we
    need a 1D bundle to shoot a single symmetry boundary condition
    $f'(0)=0$).
    Substituting into \ef{E1} yields
     \be
     \label{E5}
     D[B_0^{n-1}(n-1)m(m-1)(m-2) +
     B_0l(l-1)(l-2)] \, y^{m(n-1)+l-3}=-y+... \, .
     \ee
   One can see that equating both algebraic functions on both sides of  \ef{E5}
 is no good: this yields a unique orbit, and not a 1D bundle (see
 below). Therefore, the only possible way is to assume that
  \be
  \label{E6}
   \tex{
  m(n-1)+l-3 < 1 \LongA \frac 3n < l<1 + \frac 3n,
 }
   \ee
   and equate to zero the square bracket on the LHS in \ef{E5}.
   This gives the following cubically algebraic
 (``characteristic") equation for the exponent $l$:
 \be
 \label{E6NN}
   \tex{
 l: \quad
 l(l-1)(l-2)= B_0^{n-2} (n-1)  \frac 3n \big(\frac 3n-1 \big)\big(2- \frac
 3n\big).
  }
 \ee
It is worth mentioning that a literal balancing $m(n-1)+l-3=1$ in
\ef{E7} just means that $D=0$ in such an \ef{E4}, i.e., this orbit
also  belongs to the above 1D bundle.

After simple transformations, this is reduced to the cubic
equation
 \be
 \label{E7}
  \fbox{$
  \tex{
 l: \quad H_n(l) \equiv l(l-1)(l-2)-(n-1) \big[\frac 3n \big(\frac 3n-1 \big)\big(2- \frac
 3n\big)\big]^{\frac 2n}=0.
 }
 $}
 \ee
 Here, we are looking for positive real roots $l$ of \ef{E7}
 satisfying two conditions in
\ef{E6}.

\begin{figure}
\centering
\includegraphics[scale=0.65]{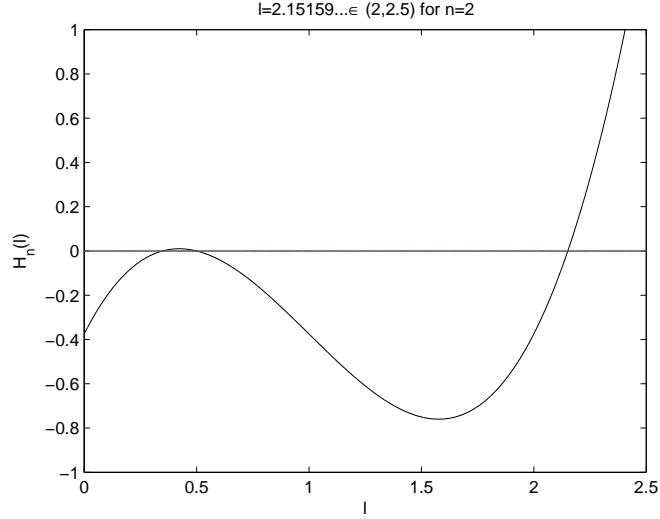}                     
 \vskip -.4cm
\caption{\rm\small The graph of $H_2(l)$ in (\ref{E7})  (for
$n=2$): $l=2.15159...$\,.}
 \label{F10}
\end{figure}

We begin our study of a proper solvability of the ``nonlinear
characteristic" equation \ef{E7} with the analytic case $n=2$, as
shown in Figure \ref{F10}, where we present the graph of the
algebraic function $H_2(l)$. Then we have
\be
\label{E8}
 \tex{
  n=2: \quad l=2.15159... \in \big(\frac 32, \frac 52
\big),
 }
  \ee
so that conditions in  \ef{E6} holds. Note that \ef{E7} shows that
there always exists an $l>2$, and, since $2> \frac 3n$ in our
parameter range, $l> \frac 3n$ always.

\begin{figure}
\centering
\includegraphics[scale=0.65]{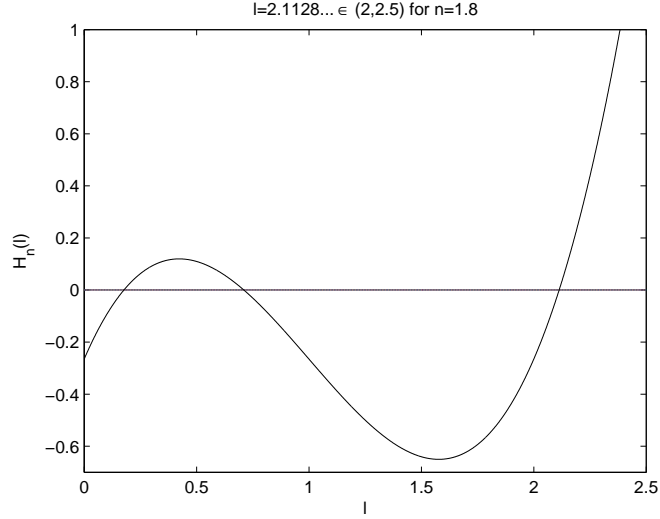}                     
 \vskip -.4cm
\caption{\rm\small  The graph of $H_n(l)$ in (\ref{E7})  for
$n=1.8$: $l=2.1128...$\,.}
 \label{F11}
\end{figure}

For $n=1.8$ (still non-oscillatory), the value of $l$ is explained
in Figure \ref{F11}, where
 \be
 \label{E9}
 n=1.8: \quad  l=2.1128...\, ,
 \ee
 also satisfying conditions in  \ef{E6}.

\begin{figure}
\centering
\includegraphics[scale=0.65]{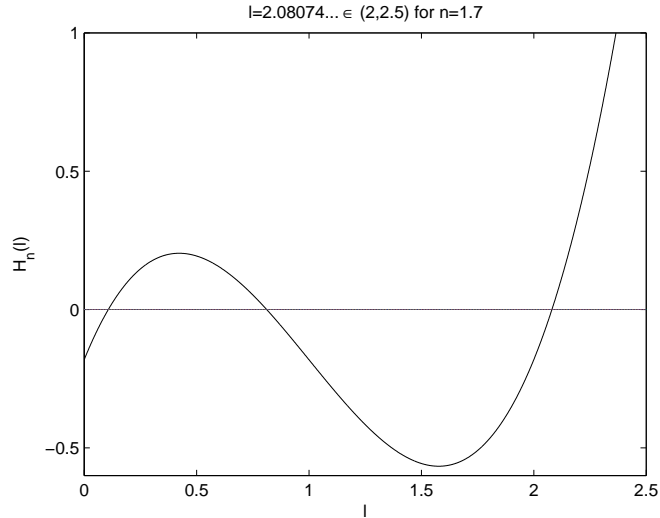}                     
 \vskip -.4cm
\caption{\rm\small The graph of $H_n(l)$ in (\ref{E7})   for
$n=1.7< n_{\rm h}$: $l=2.08074...$\,.}
 \label{F12}
\end{figure}

Finally, we have got that such a positive 1D bundle perfectly
exists also in the oscillatory range, i.e., for $n<n_{\rm h}$: see
Figure \ref{F12}, where
 \be
 \label{E10}
 n=1.7 < n_{\rm h}: \quad l=2.08074...\, .
  \ee


Concerning a proper justification of the actual existence of the
bundle \ef{E4}, \ef{E7} on a small interval $y \in [0,\d]$, on one
hand, this can be done by reducing the third-order non-autonomous
ODE  \ef{E1} to a dynamical system in $\re^4$. On the other hand,
 since the right-hand side of
  \be
  \label{KK1}
 \tex{
 f'''= - \frac {1-y}{f^{n-1}}
 }
 \ee
 is singular at $y=0$, there occurs an integrability problem.
 Indeed, a direct first integration of \ef{KK1} over a small interval $(0,y)$
 is impossible since, in the present parameter $n$-range, $f''(0)=
 +\iy$.  Therefore, one has to use another representation of the
 equation.

 To clarify this, consider, e.g., the case $n=2$. We then integrate \ef{E1} over $(0,y)$ to get
 \be
 \label{KK2}
 \tex{
 f f'' - \frac 12 \, (f')^2=-\big(y- \frac {y^2}2 \big),
  }
  \ee
  where all accompanying integrals are convergent at $y=0$. Hence,
   \be
   \label{KK3}
    \tex{
    f''= \frac 1 f \, \big[ \frac 12\, (f')^2 -\big(y- \frac {y^2}2
    \big)  \big].
    }
    \ee
   Note that,  in the class of solutions close to \ef{E4} (with $m= \frac 32$
    for $n=2$), both terms on the right-hand side of \ef{KK3} are
    equally involved into the behaviour of $f(y)$ as $y \to 0^+$.
Integrating \ef{KK3} two times, with all integrals convergent,
 since $m =\frac 3 n= \frac 32 >1$, so $f'(y)$ is locally bounded, we
arrive at a system for $\{f,f'\}$ of a standard form
 \be
 \label{KK56}
  \left(
  \begin{matrix}
  f
  \\
  f'
  \end{matrix}
  \right)
  = {\bf M}
  \left(
  \begin{matrix}
  f
  \\
  f'
  \end{matrix}
  \right),
  \ee
 with an integral operator ${\bf M}$, easily reconstructed from \ef{KK3}, being
a contraction in $C([0,\d])$ for sufficiently small $\d>0$
(meaning $C^1$ for $f$) in a properly chosen invariant
neighbourhood of each orbit from \ef{E4} for any fixed $D \in
\re$. As a ``measure" of this neighbourhood, one can use the third
term in \ef{E4}, i.e., that one appeared for $D=0$ (we have
mentioned it above).

\ssk

One can see from \ef{E7} that a proper characteristic root $l$
exists for all $n \in \big( \frac 32,3 \big)$, so that a positive
expansion \ef{E4}, \ef{E7} close to interfaces for the Cauchy
problem always exists in this whole parameter range. Moreover, one
can see that, in this parameter range,
 \be
 \label{E11}
  \tex{
  n \in \big( \frac 32, 3\big): \quad f(0)=f'(0)=(f^n f''')(0)=0,
  }
  \ee
  i.e., standard zero-high,  zero-contact angle, and zero-flux conditions, which are requirements of the FBP
  (for the CP, exhibiting, usually, a smoother behaviour at the
  interface, those are also valid), take place.

 We, thus, naturally, come to our:

\section{Final conclusions}

Thus, we have seen that the positive expansion \ef{E4}, \ef{E7}
exists for all $n$ in our range $n \in \big( \frac 32,2]$, and
also for all $n \in(2,3)$, though this range is out of our current
interests.

So, the following natural question arises: why then do we need to
take into account an oscillatory behaviour for $\frac 32 < n<
n_{\rm h}$? The answer is as follows:

\ssk

 {\bf 1.} It turns out that the shooting, via the positive 1D
bundle \ef{E4}, \ef{E7} posed at $y=y_0^-$, the single symmetry
condition at the origin $f'(0)=0$ is consistent not for all $n$'s,
and only for $n \ge n_{\rm h}$.

\ssk

 {\bf 2.} So, there exists a heteroclinic bifurcation point
$n=n_{\rm h}$, starting with which using positive expansion near
the interface becomes inconsistent, so a new oscillatory 1D bundle
is necessary, which allows to shoot a proper similarity profile
 $f(y)$.

A naive analogy in linear spectral non-self-adjoint theory is: a
spectral parameter $\l$, starting  at some critical value of a
parameter, must become complex, with a non-zero imaginary part
(meaning oscillations of eigenfunctions), in order to make
consistent algebraic systems arising.

As a consequence and as an answer of a possible (quite reasonable)
 critics concerning a definite lack of rigorous results in the
 present paper on the actual global existence of the similarity
 profiles $f(y)$ for the Cauchy problem in our range $n \in (
 \frac 32,2]$, it is worth mentioning that any kind of
 shooting\footnote{E.g., a naive direct approach using \ef{E4},
 \ef{E7}: for $D \gg 1$, we have $f'(0)>0$, and, for $D \ll -1$, we get
 $f'(0)<0$, so that, it seems, there must exist a $D^*$ in between
 such that $f'(0)=0$. But this does not work, in general, since,
 for $n < n_{\rm h}$ slightly, no such solutions were observed.
 Clearly and in addition, the above primitive analysis cannot reveal existence of
 any $n_{\rm h}$.}
 using the positive or oscillatory bundle at the interface {\sc
 must reveal this critical exponent $n_{\rm h}$}.
 Recall that $n_{\rm}$ comes from the ODE \ef{m=2.11} by using a
 ``microscopic blow-up" analysis of the ODE \ef{up2} for $f(y)$. We do not
 think that any king of a global-type study of this ODE \ef{up2}
 (without a blow-up one as $y \to y_0^-$) can reveal this critical
 heteroclinic exponent.

\ssk

{\bf 3.} Finally, bearing in mind \ef{E11} and the uniqueness of
an admissible (positive) expansion, we end up with the following
expected conclusion:
 \be
 \label{E14}
  \fbox{$
 \mbox{for $n \in [n_{\rm h},3)$, the FBP and CP have the same
 solutions,}
 $}
 \ee
 which was already announced in our earlier papers, but, unfortunately, without
 sufficiently convincing details.

\smallskip

{\bf Acknowledgements.}  The author would like to thank J.D.~Evans
 for  discussions and good questions, essentially initiated the present research, to
  be continued for the sixth-order thing film equation (the TFE--6;
 in lines with our previous papers \cite{GBl6}, where all asymptotic/numerical tools get more complicated).


\end{document}